\documentclass[12pt]{article}
\title{Beyond Dominated Convergence:
\newline
Newer Methods of Integration $\;\;\;\;\;$}
\date{}

\newtheorem{theorem}{Theorem}

\newcommand{\vt}{\vspace{5pt}\\}

\newcommand{\R}{\mathbf{R}}

\newcommand{\D}{\mathcal{D}}

\newcommand{\Pa}{\mathcal{P}}

\newcommand{\g}{\gamma}

\newcommand{\ve}{\varepsilon}

\usepackage[pdftex]{graphicx}

\author{Pat Muldowney}

\begin{document}
\maketitle
\begin{abstract}
Lebesgue's dominated convergence theorem is a crucial pillar of modern analysis, but there are certain areas of the subject where this theorem is deficient. Deeper criteria for convergence of integrals are described in this article.
\end{abstract}

\section{What Is Integration?} 
We learn calculus in secondary school: first, differentiation of functions, and later integration as the inverse or opposite of differentiation---the integral is the anti-derivative or primitive function, from which \emph{definite integrals} can be easily deduced.

In more advanced mathematics courses we learn Riemann's definition of definite integrals which enables us to integrate more functions. The Riemann method does not make use of differentiation; it is similar to the ancient method of finding areas and volumes ``by exhaustion''---that is, estimating the area or volume by dividing it up into pieces which are more easily estimated, and then taking the aggregate of the pieces. 

Specialists in mathematical analysis go on from this to study the Lebesgue method of integration. Why?
Again, one of the stock answers to this question is that the Lebesgue method enables us to integrate functions which cannot be integrated by more familiar methods such as the calculus integral and the Riemann integral. 

The Dirichlet function is sometimes mentioned. In the unit interval $[0,1]$ this function takes value one at the rational points, and zero at the other points. The Dirichlet function is not the derivative of some other function, so it cannot be integrated by the method we learn at school in calculus lessons. Also, it cannot be integrated by Riemann's method. But the Lebesgue integral of the Dirichlet function exists: the definite integral of the Dirichlet function on the unit interval is $0$.

But---so what? Apart from some specialists and experts, is there anybody else who has any real use for the Dirichlet function, and who really cares whether or not it is amenable to calculus? It cannot be pictured as lines in a graph, it does not have a straightforward formulation in polynomial, trigonometric or exponential terms. Unlike the area and volume calculations of antiquity, and unlike the calculus of Newton and Leibnitz which explained the world in mechanical terms, what {difference} does the Dirichlet function make to anyone outside the narrow and rarified world of a tiny number of people in pure mathematics?

The same can be said of many of the other arcane and exotic functions, such as the Devil's Staircase, invented during the long nineteenth century gestation of Lebesgue integration, measure theory and set theory. Such functions have counter-intuitive and challenging qualities that we can admire and wonder at. But they were described by Hermite and Poincar{\'e} as unwelcome monsters causing mayhem in the rich and fertile garden of mathematics. 

\emph{``Does anyone believe that the difference between the Lebesgue and the Riemann integrals can have physical significance, and that whether, say, an airplane would or would not fly could depend on this difference? If such were claimed, I should not care to fly in that plane.''} (Richard W.~Hamming \cite{Hamming}).

This critique is understandable, but unhistorical. By the early nineteenth century, the rich and fertile garden was on the verge of becoming a barren and dangerous wilderness---and not because of trespassing monster-functions. 

\section{Monstrous Functions}
In the tradition of Newton and Leibnitz,  Fourier series representation of functions had opened up the analysis of wave motion, crucial to an understanding of sound, light, electricity and so on. But strange and paradoxical things could happen when the integral of a function was obtained by integrating its Fourier series term by term. Certain questions could no longer be avoided. To what extent, and under what conditions, is a function identical to its corresponding Fourier series representation?
 When is the integral of a function equal to the series obtained by integrating the terms of the corresponding Fourier series? 
 
This boils down to whether a convergent series of integrable functions has integrable limit, and whether the integral of the limit is equal to the limit of the integrals  whenever the latter limit exists. 

Such issues motivated decades of investigation of the notion of integration, until a satisfactory resolution was found in the convergence theorems---uniform, monotone, and dominated convergence---of Lesbesgue integration theory. In particular the dominated convergence theorem tells us that if a sequence of integrable functions $f_j$ converges to $f$, and if the sequence satisfies $|f_j| < g$ for all $j$, where $g$ is integrable, then $f$ is integrable and $\int f_j$ converges to $\int f$ as $j \rightarrow \infty$. 

The integral here is the definite Lebesgue integral on some domain. But the theorem holds for functions which are integrable in the older and more familiar senses of Riemann, Cauchy, and Newton/Leibnitz, since, broadly speaking, functions which are integrable in the latter senses are, \emph{a fortiori}, integrable in the Lebesgue sense.

From this point onwards somebody---not necessarily expert in Lebesgue's integration---who is contentedly doing some familiar integral operations, and who encounters some issue of convergence such as term-wise integration of a Fourier series, can proceed in safety if a dominant integrable function $g$ can be found for the convergence.

This is the practical significance of Lebesgue's theory. It is a reason why ``it is safe to fly in airplanes'', so to speak. It is why the fertile garden did not turn into a barren wilderness. And the ``monster-functions'' were in reality guard dogs that played their part in protecting the garden.

\section{A Non-monstrous Function}
But is this the end of the story? Did  Lebesgue's 1901 and 1902 papers \cite{Lebesgue1, Lebesgue2} give the last word on the subject?

Here is a sequence of non-monstrous functions formed by combining some familiar polynomial with trigonometric functions. For $j=2,3,4, \ldots$,  let $f_j(x)$ $=$ $2x \sin \frac 1{x^2} - \frac 2x \cos \frac 1{x^2}$ if $\frac 1j \leq x\leq 1$ and $=$ $0$ if  $0\leq x <\frac 1j$.

An impression of function $f_j$ can be obtained from Figure 3 below, which contains the graph of  $f(x)=2x \sin \frac 1{x^2} - \frac 2x \cos \frac 1{x^2}$ for $0<x \leq 1$. Figure 2 resembles Figure 3 in the neighbourhood of $x=0$. The difference between the two is $2x \sin x^{-2}$, whose graph is in Figure 1. 

But the values of the latter function are very small in the neighbourhood of $x=0$. Figure 1 demonstrates its ``visual insignificance'', so to speak. Note that the vertical scale of Figure 1 is much more magnified than the vertical scales of Figures 2 and 3. 

Figure 4 is the graph of the primitive function or anti-derivative of $f$, which will play a big part in our discussion.

Each $f_j$ has a single discontinuity (at $x= \frac 1j$), but is  differentiable at every other point. Each $f_j$ is integrable (in the sense of Riemann and Lebesgue), and the sequence $f_j$ is convergent at each $x$ to $f(x) = 2x \sin \frac 1{x^2} - \frac 2x \cos \frac 1{x^2}$, $f(0) =0$, whose graph is Figure 3. 

{The limit function $f$ has a single discontinuity (from the right) at $x=0$; and it has a primitive function $F(x)=x^2\sin x^{-2}$ ($x>0$), $F(0)=0$ (Figure 4); so in fact $f$ {has} a definite integral $F(1)-F(0) = \sin 1$ on the domain $[0,1]$ \textbf{provided} the Newton-Leibnitz definition of the integral is used.} 

But $f$ is unbounded on $[0,1]$ and therefore is not Riemann integrable on $[0,1]$. And, though clearly non-monstrous, and understandable to a beginning calculus student, the function $f$ is \textbf{not} Lebesgue integrable.  See below for discussion of this point.

On the face of it, this example indicates a step backwards, as it were, where the old school method of Newton/Leibnitz is actually more effective than more modern methods. Leb\-esgue's theory of the integral threw up anomalies of this kind, and accordingly investig\-ation of the theory continued through the twentieth century. 

To recapitulate, Lebesgue's dominated convergence theorem can be said to be the cutting edge of modern integration theory. But it fails to capture the convergence of  sequences such as $f_j$ and $\int_0^1 f_j(x) \,dx$. The graph in Figure 3 suggests $2x^{-2}$ as a conceivable candidate for dominating function $g$ for the terms $|f_j|$, $j=1,2,3, \ldots$, at least in a neighbourhood of the critical point $x=0$. But $2x^{-2}$ is not integrable in a neighbourhood of $0$, and it seems that the dominated convergence theorem is not applicable here.

This failure must sound some alarm bells, because while many working mathematicians can get by without the Lebesgue integral, we cannot really do without convergence theorems which allow us, for instance, to perform routine operations on Fourier series; or, more generally, to safely find the integral of the limit of a sequence of integrable functions by taking the limit of the corresponding sequence of integrals.
 
The purpose of this essay is to dip into some aspects of modern integration theory in order to introduce Theorems \ref{crit1}, \ref{crit2}, and \ref{crit3} below, which are delicate enough to deal with, for instance, the convergence of the functions $f_j$ above and their integrals; while also covering the ground already captured by the convergence theorems of Lebesgue's theory.

For ease of reference, here again are  the sequence $f_j$ and related functions, $j=1,2,3, \ldots$:

\begin{eqnarray}\label{fndefs}
f(x)&=& \left\{
\begin{array}{lll}
2x \sin \frac 1{x^2} - \frac 2x \cos \frac 1{x^2} &\mbox{if} & 0<x\leq 1, \\
0&\mbox{if} & x=0,
\end{array} \right. 
\label{f}
\vt
f_j(x)&=& \left\{
\begin{array}{lll}
2x \sin \frac 1{x^2} - \frac 2x \cos \frac 1{x^2} &\mbox{if} & \frac 1j \leq x\leq 1, \\
0&\mbox{if} & 0\leq x <\frac 1j,
\end{array} \right. 
\label{fj}
\vt
F(x)&=& \left\{
\begin{array}{lll}
x^2 \sin \frac 1{x^2} &\mbox{if} & 0< x\leq 1, \\
0&\mbox{if} & x=0 ,
\end{array} \right.
\label{g}
\vt
F_j(x)&=& \left\{
\begin{array}{lll}
x^2 \sin \frac 1{x^2} &\mbox{if} & \frac 1j \leq x\leq 1, \\
0&\mbox{if} & 0 \leq x < \frac 1j.
\end{array} \right. 
\label{gj}
\end{eqnarray}

Figures 3 and 4 are, respectively, the graphs of the functions $f$ and $F$. The graphs of $f_j$ and $F_j$ are easily substituted---just insert a horizontal line segment at height $0$ from $x=0$ to $x=\frac 1j$.

The function $f$ has a single discontinuity at $x=0$, while $F$ is continuous. The reader should verify that $F$ is differentiable, including at the point $x=0$ (from the right), and that $F'(0) =0 = f(0)$. The fact that $F'(x) = f(x)$ for $x>0$ is easily verified.

This establishes that $f$, $=\lim f_j$, has a primitive or anti-derivative, and $F$ is the calculus- or Newton-indefinite integral of $f$. Also the calculus- or Newton-definite integral on $[0,1]$ is
\[
\int_0^1 f(x)\,dx = F(1) - F(0) =\sin 1 - 0 = \sin 1.
\]
For each $j$ both $f_j$ and $F_j$ have a discontinuity at $x=1/j$, but provided $x \neq 1/j$, we have $F'(x) = f_j(x)$. Thus, for each $j$, $f_j$ is Riemann and Lebesgue integrable on $[0,1]$, but not calculus- or Newton-integrable on $[0,1]$. However $f_j$ \textbf{is} calculus-integrable on $[j^{-1}, 1]$ for each $j$, and
\[
\int_0^1 f_j(x)\,dx = \int_\frac 1j ^1 f_j(x)\,dx= F(1) - F\left(\frac 1j \right) = \sin 1 - \frac 1{j^2} \sin j^2,
\]
provided we interpret $\int_0^1$ as a Riemann (or Lebesgue) integral. Thus, as $j \rightarrow \infty$,
\begin{equation}\label{convergence example}
\int_0^1 f_j(x)\,dx \rightarrow \int_0^1 f(x)\,dx
\end{equation}
provide we interpret the left hand integrals in the sense of Riemann or Lebesgue, and the right hand one as a calculus or Newton/Leibnitz integral. 

Unless we can interpret it in some other way, (\ref{convergence example}) is deficient as it stands, since we cannot ascribe the same meaning to the symbol $\int_0^1$ in the left- and right-hand terms. In fact we will establish later that (\ref{convergence example}) is valid for an adapted\footnote{That is, the Riemann-complete integral, also known as the generalized Riemann or Henstock-Kurzweil integral.} version of the Riemann integral; and that the convergence---including integrability of the limit function $f$---though unrelated to any kind of dominated convergence, satisfies a new kind of Riemann sum convergence criterion which goes beyond the Lebesgue dominated convergence theorem.









\begin{figure}
\begin{minipage}[b]{0.45\linewidth}
\centering
\includegraphics[width=\textwidth]{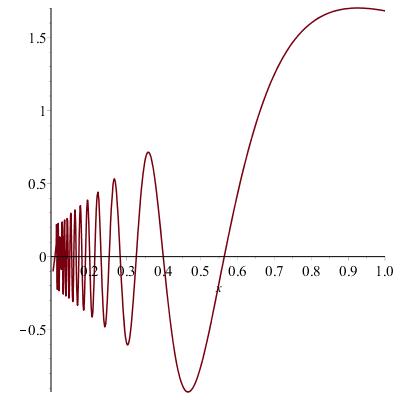}
\caption{$2x\sin x^{-2}$}
\label{fig:figure1}
\hspace{0.5cm}
\centering
\includegraphics[width=\textwidth]{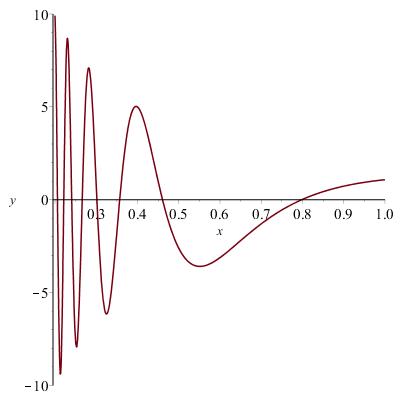}
\caption{$2x^{-1} \cos x^{-2}$}
\label{fig:figure3}

\end{minipage}
\hspace{0.5cm}
\begin{minipage}[b]{0.45\linewidth}
\centering
\includegraphics[width=\textwidth]{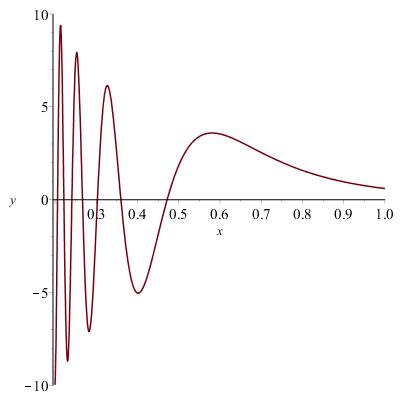}
\caption{$2x\sin x^{-2} - 2x^{-1} \cos x^{-2}$}
\label{fig:figure2}
\hspace{0.5cm}
\centering
\includegraphics[width=\textwidth]{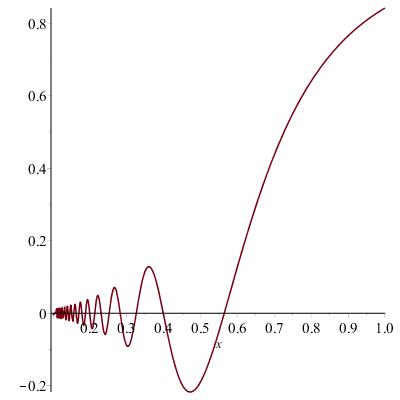}
\caption{$x^2 \sin x^{-2}$}
\label{fig:figure4}
\end{minipage}
\end{figure}

To recapitulate, for $0 \leq x \leq 1$ the function $f_j$ is bounded and continuous---except for a discontinuity at $x=j^{-1}$. Also it is differentiable except at $x=j^{-1}$. It has anti-derivative $F_j(x)$---except at $x=j^{-1}$.

By familiar standard results, $f_j$ is Riemann integrable and Lebesgue integrable on domain $[0,1]$. But, as discussed above, its limit function $f$ is \textbf{not} Lebesgue (or Riemann) integrable.

There are theorems which tell us when to expect Lebesgue integrability of the limit function of Lebesgue integrable functions. If the convergence of the functions $f_j$ to the function $f$ is uniform, monotone, or dominated by a Lebesgue integrable function, then Lebesgue integrability of the functions $f_j$ implies Lebesgue integrability of their limit function $f$, with
\[
\lim_{j\rightarrow \infty}\left(\int_0^1 f_j(x) dx\right) = 
\int_0^1 \left(\lim_{j\rightarrow \infty}f_j(x)\right) dx =
\int_0^1 F(x) dx.
\]
Inspection of the graphs indicates that convergence of functions $f_j$ is not uniform, monotone or dominated. And, even though each $f_j$ is Lebesgue integrable, the limit function $f$ is not Lebesgue integrable---as demonstrated below.

Is this a big problem for the garden of mathematics, or is it just a minor incursion by, not a monster, but an atypical creature which is easily contained? This article attempts to provide some perspective.

The names of Denjoy, Perron, Kolmogorov and others are associated with twentieth century efforts \cite{Gordon} to 
pursue the implications of problems such as the convergence of functions $f_j$ and their integrals. This article will examine the Riemann sum approach of R.~Henstock and J.~Kurzweil. 

\section{Riemann-complete Integration}
Kurzweil came to this subject through his investigations of differential equations. Henstock was interested in convergence issues in integration. Independently, each of them  focussed on careful construction of Riemann sums for integrands $f$.

Here is a broad outline of  Riemann sum construction. A partition $\Pa$ of a domain such as $[0,1]$ is a set of points $u_0<u_1<u_2< \cdots < u_n=1$. Identify $\Pa$ with the corresponding set of disjoint intervals $I_i$:
\[
\Pa = \left\{[u_0, u_1],\;]u_1,u_2],\;]u_2,u_3],\ldots , ]u_{n-1},u_n]\right\} = \{I_i\;:\;i=1,2, \ldots ,n\}.
\]
For each $I_i \in \Pa$ let $|I_i|$ denote the length $u_i - u_{i-1}$ of $I_i$. Given a function $f(x)$ defined for $x\in [0,1]$,  evaluation points $x_i$ are selected for the intervals $I_i\in\Pa$ in accordance with certain rules (such as $u_{i-1}\leq x_i <u_i$),  and then a Riemann sum for $f$ is \[
\sum_{i=1}^n f(x_i)\times |I_i|, \;\;\;\mbox{ or, more briefly, }\;\;\;(\Pa) \sum f(x)|I|.
\]
The integral $\int_0^1 f(x)\,dx$ of $f$ on the domain $[0,1]$, denoted by $\alpha$, exists if $\alpha$ satisfies a condition which is broadly of the following form.
Given $\ve>0$,  partitions $\Pa$ can be chosen such that
\begin{equation}\label{Riemann sum inequality}
\left|\alpha -  \sum_{i=1}^n f(x_i) |I_i|\right| < \ve\;\mbox{ for specified choices of }\;x_i \mbox{ and } I_i,\;1 \leq i \leq n.
\end{equation}
This inequality is reminiscent of the Riemann integral of $f$, but it is \textbf{not} the full definition that is required here. There must be some rule (sometimes called a \emph{gauge}) for selecting the elements $\{(x_i,I_i)\}$ (corresponding to the partition $\Pa =\{I_i\}$ or $\{(x_i,I_i)\}$) that can be admitted in the inequality. 

For Riemann integration, the rule is that, given $\ve>0$ there exists a constant $\delta>0$ such that, for every partition $\Pa=\{(x_i,I_i)\}$ for which $|I_i|=u_i-u_{i-1}<\delta$ and $u_{i-1}\leq x_i\leq u_i$ for each $I_i \in \Pa$, the above\footnote{If $f$ is continuous on $[0,1]$ then its Riemann integral exists there.} inequality (\ref{Riemann sum inequality}) holds. Denote such a rule by $\g$, and denote a partition $\Pa$ which satisfies an appropriate instance of the rule $\g$ by $\Pa_\g$.

An integral constructed from such a rule can be identified by notation ${^\g}\!\!\int$. 
Then the definition of the integral $\alpha =\;{^\g}\!\!\int_0^1 f(x) \,dx$ is as follows. 
There is a number $\alpha$ for which, given any $\ve>0$,  there exists a corresponding\footnote{For the ordinary Riemann integral, read ``there exists a corresponding number $\delta>0$''.} instance $\g(\ve)$ of $\g$ such that every partition $\Pa_{\g(\ve)}$ satisfies
\begin{equation}
\label{defint}
\left|\alpha - (\Pa_{\g(\ve)}) \sum f(x)|I|\right| < \ve.
\end{equation}
We will omit the $^\g$ in ${^\g}\!\!\int$, and allow the context to demonstrate which version of the integral is being discussed.

While the Riemann sum rule for ordinary Riemann integration is ``$|I|<\delta$'', the primary innovations of Kurzweil and Henstock were:
\begin{enumerate}
\item
to replace selection of intervals $\{I_i\}$ by selection of linked pairs $\{(x_i,I_i)\}$ in constructing Riemann sums,\footnote{The familiar condition $u_{i-1}\leq x_i\leq x_i$ is sometimes altered. Also the new approach often gives priority to the evaluation points $x_i$; and then it can be a more subtle and difficult task to determine the linked partitioning elements $I_i$. See \cite{The Calculus and Gauge Integrals} and \cite{MTRV}.}  and
\item
 to replace the constant $\delta$ by variable $\delta (x)>0$, where $x=x_i$ is the evaluation point in the term $f(x_i)|I_i|$ of the Riemann sum in (\ref{defint}). 
\end{enumerate} 
 To distinguish this from the Riemann integral, call it\footnote{It is also called the Henstock-Kurzweil integral, generalized Riemann integral, and gauge integral.} the \emph{Riemann-complete} integral. Clearly, every Riemann integrable function is also integrable in the Riemann-complete sense.

A Stieltjes-type definition of the integral of a function $f$ can be expressed as follows. Suppose $f(x)$ and $g(x)$ are point functions defined on the domain $[0,1]$.
The Riemann-Stieltjes integral of $f$ with respect to $g$ is got by replacing the length function $|I|$ by the increment function $g(I)=g(u_i) - g(u_{i-1})$ in the above definitions.
A standard result is that if $f$ is continuous and $g$ is monotone (or has bounded variation), then the Riemann-Stieltjes integral $\int_0^1 f(x)\,dg$ exists. If the constant $\delta>0$ in the definition is replaced by the function $\delta(x)>0$ of the Riemann-complete construction, call the resulting integral the \emph{Stieltjes-complete} integral of $f$ with respect to $g$.

To see how the Riemann-complete integral matches the calculus or Newton/Leibnitz integral, suppose a point function $f(x)$ has an anti-derivative or primitive function $F(x)$ for $0\leq x \leq 1$, so its definite integral in the Newton/Leibnitz sense is $F(1)-F(0)$. Proceeding as follows, it is easy to deduce that $f$ is Riemann-complete integrable. 

If $\Pa$ is a partition of $[0,1]$ with partition points $u_i$, $0=u_0<u_1< \cdots <u_n=1$, and if $u_{i-1} \leq x_i \leq u_i$, then
\[
\begin{array}{rll}
(\Pa)\sum f(x)|I|& = &\sum_{i=1}^n f(x_i)(u_i - u_{i-1})\vt
&=& \sum_{i=1}^n \left(f(x_i)(x_i - u_{i-1})+f(x_i)(u_i - x_i)\right), \vt
F(1) - F(0) &=& \sum_{i=1}^n F(u_i) - F(u_{i-1}) \vt
&=& \sum_{i=1}^n \left((F(x_i)-F(u_{i-1}) + (F(u_i) - F(x_i))\right)
\end{array}
\]
Let $\ve>0$ be given. Then, for each $x$, $0<x<1$, there exists a number $\delta(x)>0$ such that, for $|x-a| < \delta(x)$, 
\begin{equation}
\label{deriv}
\left|\frac{F(x) - F(a)}{x-a} - f(x)\right|< \ve.
\end{equation}
Now choose a partition $\Pa$ so that each term $f(x_i)(u_i-u_{i-1})$ satisfies 
\[
x_i - u_{i-1}< \delta(x_i),\;\;\;\;\;\;u_i -x_i < \delta(x_i).
\]
The existence of such partitions is a consequence of the Heine-Borel theorem. For such a partition (\ref{deriv}) implies
\begin{eqnarray}
\left|\left(F(x_i) - F(u_{i-1})\right) - f(x_i)(x-u_{i-1}) \right|&<&  \ve(x_i-u_{i-1}),\label{u(i-1)}\vt
\left|\left({F(u_i) - F(x_i)}\right) - f(x_i)(u_i -x_i)\right| &<& \ve(u_i-x_i);\label{ui}
\end{eqnarray}
with corresponding inequalities for $x=0$ and $x=1$. Writing $\alpha = F(1)-F(0)$ and $R=(\Pa) \sum f(x)|I|$,

\begin{eqnarray*}
\left|\alpha - R \right|&=& \left|\sum_{i=1}^n \left(F(u_i) -F(x_i)\right) +\left(F(x_i) - F(u_{i-1})\right) \;\;-  \right. \vt
&&\;\;\;\; \;\;\;\; \;\;\;\; -\;\;
\left.\sum_{i=1}^n \left(f(x_i)(x_i - u_{i-1})+f(x_i)(u_i - x_i)\right)\right| \vt
&\leq & \sum_{i=1}^n \left\{ \left|(F(x_i)-F(u_{i-1})-f(x_i)(x_i-u_{i-1})\right| \right. \;\;+ \vt
&&
\;\;\;\; \;\;\;\; \;\;\;\; +\;\;
\left.\left| (F(u_i)-F(x_i)) - f(x_i)(u_i - x_i)\right|\right\} \vt
&<& \sum_{i=1}^n \left\{ \ve(x_i - u_{i-1}) + \ve (u_i - x_i)\right\} \vt
&=&\ve \sum_{i=1}^n (x_i - u_{i-1} + u_i - x_i) \;\;\;=\;\;\;\ve.
\end{eqnarray*}
Thus the Riemann-complete integral of $f$ exists and equals the Newton/Leibnitz definite integral $F(1)-F(0)$.

In the case that $f$ is given by (\ref{f}), while  the Newton/Leibnitz and Riemann-complete integrals exist, it has been asserted above that the Lesbesgue integral does not exist.\footnote{In that case the Riemann integral of $f$ does not exist either.} 
This assertion remains to be demonstrated.

The definition of the Lebesgue integral of a function can be addressed in various equivalent ways, e.g.~\cite{Royden}. For instance, given a real-valued, measurable function $f$ defined on an arbitrary measurable space $S$, with measure $\mu$ defined on the family of measurable subsets of $S$, \cite{MulStoch} shows how to define the Lebesgue integral of  $f$ on $S$ as a Riemann-Stieltjes integral in $\R$, the set of real numbers. In fact, writing \[
g(x) = \mu\left(f^{-1}(]-\infty,x])\right),
\] 
the Lebesgue integral $\int_S f\,d\mu$ is the Riemann-Stieltjes integral $\int_{-\infty}^\infty f\,dg$. If $S \subseteq \R$ and $\mu$ is Lebesgue measure in $\R$, and if the Lebesgue integral of $f$ exists, then the Riemann-complete integral of $f$ exists and the two integrals are equal. Every Lebesgue integrable function is integrable in the Riemann-complete sense (see \cite{MTRV}).

A key point is that the Lebesgue integral is an absolute integral, while the Riemann-complete is non-absolute. Writing $f_+(x) = f(x)$ if $f(x) \geq 0$, with $f_+(x) =0$ otherwise, and $f_-(x) =f(x) - f_+(x)$, absolute integrability implies that $f$ is Lebesgue  integrable if and only if both $f_+$ and $f_-$ are Lebesgue integrable.\footnote{This restriction does not apply to the Riemann-complete integral of $f$, which does \textbf{not} require the Riemann-complete integrability of $f_+$ and $f_-$.} We use this point to demonstrate that the function $f$ defined by (\ref{f}) is not Lebesgue integrable.

With that  in mind, Figure 3 provides an indication of how the function $f$ defined by (\ref{f}) 
\begin{itemize}
\item
fails to be Lebesgue integrable, while
\item
its Riemann-complete integral exists.
\end{itemize}
In fact Figure 3 shows that, in neighbourhoods of $x=0$, the graph of $f$ oscillates increasingly rapidly, in positive ($f_+$, above the $x$-axis) loops and negative ($f_-$, below the $x$-axis) loops whose amplitude (or height/depth) increases without limit as $x \rightarrow 0$. This creates the suspicion, or expectation, that the sum of areas of the positive loops diverges to $+\infty$, while the sum of areas of the negative loops diverges to $-\infty$. 

But if, instead of treating positive and negative loops separately, we add up their areas  in their natural sequence, then positive and negative areas will tend to cancel each other out, and the resulting sequence of net values may converge.\footnote{For example, with $b_i = (-1)^i i^{-1}$, the series $\sum_{i=1}^\infty b_i$ converges, but the series consisting of only the positive terms (or only negative terms) diverges; so $\sum_{i=1}^\infty |b_i|$ diverges.}
The latter is what happens in the Riemann sum construction of the Riemann-complete integral of $f$. 

The following discussion seeks to add substance to these speculations.
In any interval, not including zero, but with small values of $x$, Figure 1 shows that the contribution from the term $2x\sin x^{-2}$  to the area under the graph of $f$ is vanishingly small in neighbourhoods of $x=0$, while the corresponding contribution from the term $2x^{-1}\cos x^{-2}$ in $f$ is relatively large. Therefore, disregarding the term $2x\sin x^{-2}$, the zeros of (\ref{f}) can, for present purposes,  be estimated approximately as
\[
x = \sqrt{ \frac 2{(2n+1)\pi}} \mbox{ as } x \rightarrow 0 \;\; (\mbox{or integer }n\rightarrow \infty ).
\]
Accordingly we may estimate that, for large, even values of $n$,
\[
\int_{\sqrt{ \frac 2{(2n+3)\pi}}}^{\sqrt{ \frac 2{(2n+1)\pi}}} f_+(x)\,dx\;\;
\mbox{ is approximately }\;\;
\frac {2}\pi \left(\frac 1{2n+1}+\frac 1{2n+3}\right)
\] while for large and odd values of $n$
\[\int_{\sqrt{ \frac 2{(2n+3)\pi}}}^{\sqrt{ \frac 2{(2n+1)\pi}}} f_-(x)\,dx\;\;
\mbox{ is approximately }\;\;
\frac {2}\pi \left(\frac 1{2n+1}+\frac 1{2n+3}\right)
.
\]
Writing
\[
a_n=\frac {2}\pi \left(\frac 1{2n+1}+\frac 1{2n+3}\right),
\]
each of the two series
\[
a_2+a_4+a_6 + \cdots,\;\;\;\;\;\;a_1+a_3+a_5 + \cdots
\]
diverges, so it is clear that $f$ is not Lebesgue integrable in $[0,1]$.
But it is easy to see that the series
\[
-a_1+a_2-a_3+a_4 - \cdots
\]
is non-absolutely
\index{integrability!non-absolute}%
convergent, even if we did not already know, from existence of the primitive function $F(x)$ for $0\leq x \leq 1$, that $f$ is Riemann-complete integrable.

This is because  the Riemann-complete convergence is obtained from the cancellation effects produced by successively summing the positive and negative parts in their natural sequence.

We can ensure this by choosing $\delta(x)$ as follows.
 When $x$ lies between adjacent roots $\sqrt{ \frac 2{(2n+1)\pi}}$ and $\sqrt{ \frac 2{(2n+3)\pi}}\;$ 
let
\[
\delta(x) < \min \left\{ x - \sqrt{ \frac 2{(2n+3)\pi}},\;\;\; \sqrt{ \frac 2{(2n+1)\pi}} -x \right\};
\]
and if $x $ is one of the roots $\sqrt{ \frac 2{(2n+1)\pi}}$, take
\[
\delta(x)< \min\left\{\, \sqrt{ \frac 2{(2n+1)\pi}},\;\;\;\sqrt{ \frac 2{(2n+1)\pi}} - \sqrt{ \frac 2{(2n+3)\pi}}\,\,\,\right\};
\]
and when $x=0$ let $\delta(0)>0$ be arbitrary.
Any partition corresponding to this definition of $\delta(x)$ ($0\leq x \leq 1$) will contain a term with $f(0)=0$, and the terms for non-zero $x$ will each contain an arbitrarily close estimate of the area of the corresponding positive or negative loop in Figure 3. This provides the required cancellation and convergence of Riemann sums, since the alternating loops are monotone decreasing in area as $x$ approaches $0$.

In the case of the Lebesgue integral this can\-c\-ellation effect is absent, and  convergence fails.

This establishes that, just as there are Lebesgue integrable functions that are not Riemann integrable, there are Riemann-complete integrable functions that are not Lebesgue integrable.

\section{Convergence Criteria}
Anybody experienced in the theory of integration will be aware that most of the preceding discussion covers fairly well-worn ground which has already been worked through in many excellent publications, such as \cite{Bartle}.

But at the outset of this article it was stated that, while Lebesgue's dominated convergence theorem is a crucial pillar of modern analysis, there are certain areas of the subject where this theorem is deficient. The sequence $\{f_j\}$ of (\ref{fj}) demonstrates that the dominated convergence theorem provides no illumination in this particular instance of converging non-absolute integrals. This section  addresses the deficit.

There are convergence conditions and criteria which encompass and surpass the dominated convergence, monotone and uniform convergence theorems of standard integration theory. These are the convergence criteria of Theorems \ref{crit1}, \ref{crit2}, and \ref{crit3} below. They are valid for Riemann-complete integrals (which include integrals of the Newton/Leibnitz, Riemann, and Lebesgue kinds). Measurability of the integrand functions is not assumed.

\begin{theorem}\label{crit1}
Suppose $f_j$ is integrable on $[0,1]$ and $f_j(x)$ converges to $g(x)$ for $x \in [0,1]$. 
Suppose, given arbitrary $\ve>0$, there exist a number $\alpha_1$ and, for $x \in [0,1]$, a gauge $\delta(x)$, and integers $p=p(x)$ depending on  $\ve$, so that, for every partition $\{I_i\}$ of $[0,1]$ with linked elements $\{(x_i,I_i)\}$ satisfying $|I_i|<\delta (x_i)$ ($i=1, \ldots ,n$), the condition
\[
\left| \alpha_1 -\sum_{i=1}^n f_{j(x_i)}(x_i)\times |I_i|\right|\;\; <\;\; \ve
\]
holds for all choices of $j=j(x_i) > p(x_i) $ ($i=1, \ldots , n$) in the Riemann sum. Then  the limit function $g(x)$ is integrable on $[0,1]$, with $\int_0^1 g(x)\,dx = \alpha_1$. 
\end{theorem}

\begin{theorem}\label{crit2}
Suppose $f_j$ is integrable on $[0,1]$ and $f_j(x)$ converges to $g(x)$ for $x \in [0,1]$.
Suppose, given arbitrary $\ve>0$, there exist a number $\alpha_2$ and a positive integer $q = q(\ve)$ depending only on $\ve$,  so that, for every partition $\{I_i\}$ of $[0,1]$ with linked elements $\{(x_i,I_i)\}$ satisfying $|I_i|<\delta (x_i)$ ($i=1, \ldots ,n$), the condition
\[
\left| \alpha_2 -\sum_{i=1}^n f_j(x_i)\times |I_i|\right|\;\; <\;\; \ve
\]
holds for every choice of $j > q(\ve) $ with $j$ constant for each term of the Riemann sum. 
Then $\int_0^1 f_j(x)dx$ converges as $j\rightarrow \infty$.
\end{theorem}
Theorems \ref{crit1} and \ref{crit2} can be expressed in converse form (see \cite{MTRV}), so they are criteria for their  respective conclusions.
Apart from the opening sentence of each they apply independently of each other; in the sense that either one of them may hold for particular integrands while the other one does not hold.

\begin{theorem}\label{crit3}
If both of Theorems \ref{crit1} and \ref{crit2} hold (so both of  $\int_0^1 g$ and $\lim_{j \rightarrow \infty}\int_0^1 f_j$ exist), then
\[
 \int_0^1 g(x)\,dx=\lim_{j \rightarrow \infty}\int_0^1 f_j(x)\,dx
 \]
if and only if $\alpha_1=\alpha_2$.
\end{theorem}
For anybody more familiar with the classical integration theorems on passage to a limit, these theorems or criteria may appear  somewhat indigestible  at first sight.

Their starting point is that a convergent sequence of  functions $f_j$ is given. These function are assumed to be Riemann-complete integrable, which is a weaker assumption than Lebesgue integrability. There is no assumption of properties like continuity or measurability. 

To answer questions about convergence of the corresponding sequence of Riemann-complete integrals, and about the Riemann-complete integrability of the limit function, from previous experience of integration we might be led to expect some condition, not about Riemann sums, but only about the  functions $f_j$---such as monotonicity, domination by a fixed integrable function $g$; or the like.
But nothing like this appears in the above convergence criteria. Instead we have various statements about Riemann sums. 

However, setting aside for a moment the conception of integral as primitive function, or anti-derivative, the original and more durable meaning of integral involves slicing up (partitioning), followed by summation, followed by taking a limit of the sums.

From this perspective, it may be less of a surprise that Riemann sums appear in the formulation of conditions for limits of integrals, since Riemann sums are central to the concept of integral.

Consider Theorem
\ref{crit1}.
Given integrability of the terms $f_j$ in the sequence, this theorem  addresses the integrability of the limit function $f$, which, essentially, involves the question of convergence of Riemann sums $\sum f(x) |I|$. 

To make an initial stab at this question, we might consider sequences of Riemann sums $\sum f_j(x) |I|$, $j=1,2,3, \ldots$.  We know that, for each $x$, the sequence of values $f_j(x)$ converges to $f(x)$ as $j\rightarrow \infty$. We also know that, for each $j$, Riemann sums of the form $\sum f_j(x)|I|$ converge to the integral of $f_j$. Can we somehow put these two facts together to deduce, as an immediate consequence,  convergence of Riemann sums $\sum f(x)|I|$ to the integral of $f$?

Of course, we know that the answer to this is \textbf{no}. The answer is \textbf{yes} if the terms $f_j$ satisfy some conditions such as $|f_j| <g$ where $g$ is integrable. But if we want a condition expressed in the form of a condition on Riemann sums, clearly something more delicate than convergence of $\sum f_j(x)|I|$ is required.

For instance, with $\ve>0$ given, the condition we need is \textbf{not} that,  for all
$j$ greater than some $j_0 = j_0(\ve)$, every Riemann sum $\sum f_j(x) |I|$ will be contained within some ball $B(\ve)$ of the form $]\beta-\ve, \beta+\ve[$. 

This is too crude for our purpose. All it says is that $f_j$ is integrable---which we already know. The convergence of $f_j(x)$ to $f(x)$ may be very fast at some points $x$, and very slow at other points $x$. This behaviour is provided for in Theorem \ref{crit1}, by choosing, \textbf{not} $j_0 = j_0(\ve)$, but $j_0 = j_0(\ve, x)$ different for each $x$.

This formulation is sufficient, and necessary, for integrability of the limit function $f$. Once this point is established, the criteria of Theorems \ref{crit2} and \ref{crit3} are fairly obvious, and less subtle.
But are these conditions of Theorems \ref{crit1}, \ref{crit2} and \ref{crit3} ``workable'' in the way that the dominated convergence condition $|f_j| < g$ is? 

After all, Riemann sums are fine for defining the meaning of the integral of a function. But when we actually want to find the value of an integral we do not typically work with Riemann sums. Instead we revert to the integral as primitive or anti-derivative, using the substitution method or integration by parts. Or we use some less direct method, such as solving a related differential equation; or a myriad of other\footnote{Which is \textbf{not} to say that Riemann sums  are ``merely'' a device of fundamental theory, and nothing else. Versions of them have had other uses; such as the ancient techniques of quadrature; or computer programs for estimating numerical values of an integral. Simpson's rule is another example.} \emph{ad hoc} methods.

To respond to such questions, and to demonstrate that Riemann sums \textbf{can} actually be of use here, we can as an example take the sequence $f_j$ defined in (\ref{fj}). Remember, for each $j$ the function $f_j$ is Riemann integrable and Lebesgue integrable, but not Newton/Leibnitz integrable, and  their limit function $f$ is Newton/Leibnitz integrable  but not Riemann or Lebesgue integrable. 
For each $j$, $f_j$ is Riemann-complete integrable.\footnote{As is $f$, from earlier discussion. But for present purposes we wish to \textbf{deduce} Riemann-complete integrability of $f$ from Theorem \ref{crit1}.} 
This discussion of the convergence criteria of Theorems \ref{crit1}, \ref{crit2} and \ref{crit3} is set in the context of Riemann-complete integrability.

The subject of the first criterion is the (Riemann-complete) integrability of the limit function $f$, and it is established by examining Riemann sums of the form
\[
\sum_{i=1}^n f_{j(x_i)}(x)\times |I_i|,\;\;\;\mbox{ or }\;\;
\sum_{i=1}^n f_{j(x_i)}(x)(u_i-u_{i-1}).
\]
We already know, by various means, including a direct investigation of the Riemann sums $\sum f(x) |I|$, that $f(x)$ is (Riemann-complete) integrable on $[0,1]$. 

The function $f(x)$ is the limit of functions $f_j(x)$. 
Is it possible to  confirm further the integrability of $f$ by direct examination, not just of $\sum f(x) |I|$, but of Riemann sums $\sum f_{j(x)}(x)|I|$ involving functions $f_j$ instead of $f$, where the factor $f_j$ in the sum has variable index $j=j(x)$, depending on the element $x$ of the division $\D=\{(x,I)\}$ used to construct the Riemann sum? 

This is the essence of Theorem \ref{crit1}.
And according to Theorem \ref{crit1} the answer to this question should be yes. Given $\ve>0$, and with a suitable gauge $\delta(x)$, provided factors $ f_{j(x)}(x)$ are chosen appropriately we should be able to demonstrate that the value of each corresponding Riemann sum $\sum f_{j(x)}(x)|I|$ will lie within some ball $B$ of radius $\ve$ where $\ve$ is arbitrarily small.


Since we are already convinced of the integrability of $f$ in this case, what we are really trying to do here is to get a sense of the behaviour of sums $\sum f_{j(x)}(x)|I|$. So, given the integrability of $f$, write $\alpha_1 = \int_{]0,1]} f(x)dx$ and choose a gauge $\delta$ so that, for every $\delta$-fine partition $\{(x_i,I_i)\}$ of $]0,1]$,
\[
\left| \alpha_1- \sum_{i=1}^n f(x_i)|I_i| \right| < \ve,\;\;\;\mbox{ or }\;\; \sum _{i=1}^n f(x_i)|I_i| \in B(\alpha_1,\ve),
\]
the ball with centre $\alpha_1$ and radius $\ve$.

Consider any one of these Riemann sums $\sum_{i=1}^n f(x_i)|I_i|$, corresponding to a particular $\delta$-fine partition with $\{(x_i,I_i)\;:\;i=1,\ldots , n\}$. For each $x$  choose 
\[
r(x)\geq \frac 1{x},\;\;\mbox{ so }\;\;r(x_i)\geq \frac 1{x_i}\;\mbox{ for each }\;i.
\]
Then, by definition of $f_j$, if $j=j(x) \geq r(x)$, 
\[
f_j(x) = f_{j(x)}(x) = f(x),
\]
so, for all choices $j(x_i) \geq r(x_i) = r(x_i, \ve)$,
\[
\sum_{i=1}^n f_{j(x_i)}(x_i)|I_i| =\sum_{i=1}^n f(x_i)|I_i| \in B(\alpha_1,\ve),
\]
as required by Theorem \ref{crit1}.

In general, the convenient equation $f_j(x) = f_{j(x)}(x) = f(x)$ cannot be appealed to. But if, with suitable choices of $j=j(x)$, the differences \[f_j(x)-f(x) = f_{j(x)}(x) - f(x),\] can make sufficiently small contributions to the Riemann sum, then it may be plausible that
\[
 \sum_{i=1}^n f_{j(x)}(x_i)|I_i| \in B'(\ve)\;\;\;\mbox{ implies }\;\;\sum_{i=1}^n f(x_i)|I_i| \in B''(\ve)
 \]
 so $f$ is integrable. This is the intuitive content of Theorem \ref{crit1}.

Now to Theorem \ref{crit2}. The preceding remarks are concerned with the integrability of $\lim_{j \rightarrow \infty }f_j$. The fundamental assumption is that each function $f_j$ in the sequence $\{f_j\}$ is integrable. In the case of our example (\ref{fj}) the anti-derivatives (\ref{gj}) are the sequence $\{F_j\}$, giving a sequence of integrals
 \[
 \int_{]0,1]} f_j(x) dx = F_j(1) - F_j(j^{-1}) = \sin 1  - F_j(j^{-1}),
 \]
which can be denoted by $\beta_j$. Note that continuity of $F$ implies $\beta_j \rightarrow 0$ as $j \rightarrow \infty$. In this case it is already clear that the sequence of integrals $ \int_{]0,1]} f_j(x) dx$ converges as $ j \rightarrow \infty$, the limit being (in this case) $\sin 1$; what we want is confirmation, including intuitive confirmation, that Theorem \ref{crit2} actually works.

The convergence of a sequence of integrals is the subject of Theorem \ref{crit2}, and it is again expressed in terms of Riemann sums. The criterion implies that, with arbitrarily small $\ve>0$ given, there is a ball $B(\alpha_2,\ve)$ with centre $\alpha_2$ and radius $\ve$, and a corresponding integer $q$ depending only on $\ve$, so that for each $j\geq q=q(\ve)$,  a gauge $\delta_j(x)>0$ can be found such that for every $\delta_j$-fine partition of $[0,1]$ the corresponding Riemann sum $\sum_{i=1}^n f_j(x_i) |I_i|$ is contained in  $B(\alpha_2,\ve)$; so
 for every $\delta_j$-fine $\{(x_i,I_i)\}$, 
\[
\left|\alpha_2 -\sum_{i=1}^n f_j(x_i) |I_i|  \right| < \ve
\]
whenever $j\geq q$.
Unlike Theorem \ref{crit1}, here $j$ is the same for each term of any particular Riemann sum. 

Again, this is easy to demonstrate  since we already know in this case that the integrals $ \int_{]0,1]} f_j(x) dx$ converge to the value $\sin 1$ as $ j \rightarrow \infty$. Just take 
\[
\alpha_2 = \sin 1 = \lim_{j \rightarrow \infty} g(1) - g(j^{-1}) = \int_{]0,1]} f_j(x) dx,
\]
and choose $q$ so that $j \geq q$ implies 
\[
\left| \alpha_2 - \int_{]0,1]} f_j(x) dx \right| < \ve.
\]
For each $j \geq q$ choose a gauge $\delta_j(x)$ ($0\leq x \leq 1$) so that, for any $\delta_j$-fine partition $\{(x_i,I_i)\}$ of $]0,1]$,
\[
\left| \int_{]0,1]} f_j(x) dx - \sum_{i=1}^n f_j(x_i) |I_i|\right| < \ve.
\]
Then, by the triangle inequality,
\[
\left| \alpha_2 - \sum_{i=1}^n f_j(x_i) |I_i| \right| < 2\ve,\;\;\;\mbox{ or }\;\;\sum_{i=1}^n f_j(x_i) |I_i| \in  B(\alpha_2, 2\ve)
\]
for all $j \geq q = q(\ve)$ and all $\delta_j$-fine partitions of $]0,1]$.
In other words, the criterion of Theorem \ref{crit2} confirms the convergence of the sequence of integrals \[\left\{\int_{]0,1]} f_j(x) dx\right\};\]
and this demonstration illustrates the intuitive content of Theorem \ref{crit2}.

Finally, the question arises whether the integral of the limit \[
\int_{]0,1]}\lim_{j\rightarrow \infty} f_j(x) dx
\] equals the limit of the integrals 
\[
\lim_{j\rightarrow \infty}\int_{]0,1]} f_j(x) dx.
\] 
For the sequence $f_j$ of (\ref{fj}), we already know by direct evaluation that these two quantities have the same value, namely $\sin 1$.
This agrees with the criterion of Theorem \ref{crit3}, which requires that $\alpha_1$ and $\alpha_2$ have the same value. In this case
\[
\alpha_1 = \sin 1=\alpha_2;
\]
so the intuitive content of Theorem \ref{crit3} is clear in the context of this example.

\section{Conclusion}
So, does it really matter whether aviation designers work out their aerodynamic equations using old-fashioned Riemann integrals or the latest fancy Lebesgue integrals?

Probably not much. But it matters a lot 
if the value 22/7 for $\pi$ were hard-wired into every computer in the world. Or if the wrong value for elasticity of O-rings at freezing temperature was used in space shuttle design.
And it certainly matters whether our aviation designers make tricky, unjustifiable calculations involving, for instance, term-by-term integration of Fourier series.

It is thanks to the intellectual diligence of the nineteenth century, not to mention its monster-functions, that we have the dominated convergence theorem to keep the garden of mathematics safe and fertile---and, indeed, to keep airplanes flying safely. 

But do we really need anything more than the dominated convergence theorem for absolutely convergent integrals? Why bring up the convergence criteria of Theorems \ref{crit1}, \ref{crit2}, and \ref{crit3}? Is the sequence $\{f_j\}$ described in (\ref{fj}) above just an exceptional one-off, or is it representative of something more significant? If the latter, where are all these non-absolute integrals?

In fact they are very widespread. Modern stochastic calculus \cite{Oksendal, Ross} is based on integrals for which absolute convergence fails, but which may converge weakly or, in some cases, non-absolutely. These are described in \cite{MTRV, MulStoch}.

A very significant formulation of quantum mechanics  is in terms of path integrals \cite{feynman2} which also fail to converge absolutely. Famously, the dominated convergence theorem does not work for these integrals, and, as described in \cite{MTRV}, the non-absolute convergence criteria must be invoked.

\emph{``Does anyone believe ...~I would not care to fly in that plane.''} A healthy scepticism must be welcomed. But what is certain is that, while integration is central to mathematical analysis, there are no certain and definite ways of tackling any problem of integration, and even a beginning student has to exercise imagination and ingenuity. From
the ancient methods of quadrature, to the methods of Newton/Leibnitz, Cauchy, Riemann, Lebesgue, Denjoy, Perron, Kolmogorov, Kurzweil, or Henstock, 
it is unwise to disregard any resource or insight that can be called upon.

\vspace{10pt}

\noindent
\emph{Pat Muldowney, Lisselty, Dunmore East, Co.~Waterford}

\noindent
\emph{pat.muldowney@gmail.com}

\end{document}